\documentclass[12pt]{article}
\usepackage[all]{xy}
\title{On t-extensions of abelian groups}
\author{H.Sahleh\\Department of Mathematics, University of Guilan, P.O.BOX 1914\\Rasht-Iran\\
e-mail: sahleh@guilan.ac.ir\\
A.A. Alijani\\ Department of Mathematics, University of Guilan\\e-mail: taleshalijan@phd.guilan.ac.ir }
\begin{document}
\date{}
\maketitle
\newcommand{\stk}[1]{\stackrel{#1}{\longrightarrow}}
\begin{abstract}
Let $\Re$ be the category of all discrete abelian groups and $\pounds$ the category of all locally compact abelian (LCA) groups. For a group $G\in \pounds$, the maximal torsion subgroup of $G$ is denoted by $tG$. A short exact sequence $0\to A\stackrel{\phi}{\to} B\stackrel{\psi}{\to}C\to 0$ in $\Re$ is said to be a t-extension if $0\to tA\stackrel{\phi}{\to} tB\stackrel{\psi}{\to}tC\to 0$ is a short exact sequence. We show that the set of all t-extensions of $A$ by $C$ is a subgroup of $Ext(C,A)$ which contains $Pext(C,A)$ for discrete abelian groups $A$ and $C$. We establish conditions under which the t-extensions split and determine those groups in $\Re$ which are t-injective or t-projective in $\Re$. Finally we determine the compact groups $G$ in $\pounds$ such that every pure extension $0\to G\to Y\to C\to 0$ splits where $C$ is a compact connected group in $\pounds$. \\
\end{abstract}

\centerline {Introduction}

Throughout, all groups are Hausdorff topological abelian groups and will be written additively. Let $\pounds$ denote the category of locally compact abelian (LCA) groups with continuous homomorphisms as morphisms and $\Re$ the category of discrete abelian groups. The Pontrjagin dual and the maximal torsion subgroup of a group $G\in\pounds$ are denoted by $\hat{G}$ and $tG$,respectively.
A morphism is called proper if it is open onto its image and a short exact sequence $0\to A\stackrel{\phi}{\to} B\stackrel{\psi}{\to}C\to 0$ in $\pounds$ is said to be proper exact if $\phi$ and $\psi$ are proper morphisms. In this case the sequence is called an extension of $A$ by $C$ ( in $\pounds $ ). Following [4], we let $Ext(C,A)$ denote the (discrete) group of extensions of $A$ by $C$. A subgroup $H$ of a group $C$ is called pure if $nH=H\cap nC$ for all positive integers $n$. An extension $0 \to A \stackrel{\phi}{\to} B \stackrel{\psi}{\to} C \to 0$ is called a pure extension if $\phi(A)$ is pure in $B$.  The elements represented by pure extensions of $A$ by $C$ form a subgroup of $Ext(C,A)$ which is denoted by $Pext(C,A)$. In section 1 and 2, all groups are discrete abelian groups. An extension $0\to A\stk{\phi} B\stk{\psi} C\to 0$ in $\Re$ will be called  a t-extension if $0\to tA\stk{\phi} tB\stk{\psi} tC\to 0$ is an extension. Let $Ext_{t}(C,A)$ denote the set of all elements in $Ext(C,A)$ represented by t-extensions. In section 1, we show that $Ext_{t}(C,A)$ is a subgroup of $Ext(C,A)$ which contains $Pext(C,A)$ (see Theorem 1.6 and Lemma 1.7). In section 2, we establish some results on splitting of t-extensions (see Lemma 2.1, Theorem 2.11 and Theorem 2.13).
Assume that $\Im$ is any subcategory of $\pounds$. The section four is a part of an investigation which answers the following question:\\
Under what conditions on $G\in\pounds$, $Ext(X,G)=0$ or $Pext(X,G)=0$ for all $X\in\Im$?. In [2], [3], [4], [5], [7] and [8] the question is answered in some subcategories of $\pounds$ such as the category of divisible locally compact abelian groups. In [5, Corollary 3.4], Fulp and Griffith proved that a compact group $G$ satisfies $Ext(C,G)=0$ for all compact connected groups $C$ if and only if $G\cong (R/Z)^{\sigma}$ where $\sigma$ is a cardinal. It may happen that $Ext(X,G)\neq 0$ but $Pext(X,G)=0$. For example, $Ext(Z_{n},Z)\neq 0$ but $Pext(Z_{n},Z)=0$ where $Z$ is the group of integers and $Z_{n}$ is the cyclic group of order $n$. In this paper, we show that a compact group $G$ satisfies $Pext(C,G)=0$ for all compact connected groups $C$ if and only if $G\cong (R/Z)^{\sigma}\bigoplus H$ where $H$ is a compact totally disconnected group (see Theorem 3.2). For the characterization of compact groups $G$ which $Pext(C,G)=0$ for all compact connected groups $C$, we need to show that $Pext(X,A)=0$ for a discrete torsion group $X$ and a discrete torsion-free group $A$ (see Corollary 2.2). As a result,[5, Corollary 3.4] is a consequence of Theorem 3.2 (see Remark 3.3).

The additive topological group of real numbers is denoted by $R$ and $Q$ is the group of rationals with the discrete topology. We denote by $G_{0}$, the identity component of a group $G\in \pounds$. For more on locally compact abelian groups, see [6].\\

\section{t-extension}

 In this section, we define the concept of a t-extension of $A$ by $C$. We show that the set of all t-extensions of $A$ by $C$ form a subgroup of $Ext(C,A)$ which contains $Pext(C,A)$.\\

{\bf Definition 1.1}. {\it An extension $0\to A\stk{\phi} B\stk{\psi} C\to 0$ is called a t-extension if $0\to tA\stk{\phi} tB \stk{\psi} tC\to 0$ is an extension.}\\

{\bf Lemma 1.2}. {\it A pushout or a pullback of a t-extension is a t-extension.}\\

Proof. Suppose $0\to A\stk{\phi} B\stk{\psi} C\to 0$ is a t-extension and

\[
\xymatrix{
0 \ar[r] & A\ar^{\phi}[r] \ar^{\mu}[d] & B \ar^{\psi}[r] \ar^{}[d] & C\ar[r] \ar^{1_{C}}[d]
& 0 \\
0 \ar[r] &A' \ar^{\phi'}[r] &   (A'\bigoplus B)/ H \ar^{\psi'}[r] & C \ar[r] & 0
}
\]

is a standard pushout diagram (see [1]). Then $$H=\{(\mu(a),-\phi(a)),a\in A\}$$ and $$\phi':a'\longmapsto (a',0)+H,~~~~~\psi':(a',b)+H\longmapsto \psi(b)$$ We show that $0\to tA'\stk{\phi'} t((A'\bigoplus B)/H)\stk{\psi'} tC\to 0$ is exact. First, we show that $\psi':t((A'\bigoplus B)/H)\to tC$ is surjective. Let $c\in tC$. Since $0\to tA\stk{\phi}tB\stk{\psi}tC\to 0$ is exact, so there exists $b\in tB$ such that $\psi(b)=c$. Clearly, $(0,b)+H\in t((A'\bigoplus B)/H)$. On the other hand, $\psi'((0,b)+H)=\psi(b)=c$. Hence $\psi'$ is surjective. Now, we show that $Ker\psi'\mid_{X}\subseteq Im\phi'\mid_{tA'}$ where $X=t((A'\bigoplus B)/H)$. Let $(a',b)+H\in X$ and $\psi'((a',b)+H)=0$. So, $\psi(b)=0$. Hence, there exists $a\in A$ such that $\phi(a)=-b$. On the other hand, there exists a positive integer $n$ such that $(na',nb)\in H$. So, there exists $a_{1}\in A$ such that $\mu(a_{1})=na'$ and $-\phi(a_{1})=nb$. Now, we have $$\phi(a_{1}-na)=\phi(a_{1})-n\phi(a)=0$$ So $a_{1}=na$ and $n(a'-\mu(a))=0$. It follows that $a'-\mu(a)\in tA'$ and $\phi'(a'-\mu(a))=(a'-\mu(a),0)+H=(a',b)+H$ (since $(a'-\mu(a),0)-(a',b)=(-\mu(a),-b)=(\mu(-a),-\phi(-a))\in H$). Now, suppose that
\[
\xymatrix{
0 \ar[r] & A\ar^{\phi'}[r] \ar^{1_{A}}[d] & B' \ar^{\psi'}[r] \ar^{}[d] & C'\ar[r] \ar^{\gamma}[d]
& 0 \\
0 \ar[r] &A \ar^{\phi}[r] &   B\ar^{\psi}[r] & C \ar[r] & 0
}
\]

is a standard pullback diagram. Then $$B'=\{(b,c');\psi(b)=\gamma(c')\}$$ and $$\phi':a\longmapsto (\phi(a),0),~~~~~~\psi':(b,c')\longmapsto c'$$ We show that $0\to tA\stk{\phi'}tB'\stk{\psi'}tC'\to 0$ is exact. Let $c'\in tC'$. Then, there exists a positive integer $n$ such that $nc'=0$. Since $\psi$ is surjective, $\psi(b)=\gamma(c')$ for some $b\in B$ . Now, $n\psi(b)=\gamma(nc')=0$. Hence, $\psi(b)\in tC$. Since $0\to tA\stk{\phi}tB\stk{\psi}tC\to 0$ is exact, so $\psi(b_{1})=\psi(b)$ for some $b_{1}\in tB$. Hence, $(b_{1},c')\in tB'$ and $\psi'(b_{1},c')=c'$. Therefore, $\psi':tB'\to tC'$ is surjective. Now, suppose that $(b,c')\in tB'$ and $\psi'(b,c')=0$. Then $c'=0$ and $nb=0$ for some positive integer $n$. So $b\in tB$. Since $\psi(b)=\gamma(c')=0$ and $0\to tA\stk{\phi}tB\stk{\psi}tC\to 0$ is exact, there exists $a\in tA$ such that $\phi(a)=b$. Now, we have $$\phi'(a)=(\phi(a),0)=(b,0)=(b,c')$$ It follows that $Ker\psi'\mid_{tB'}\subseteq Im\phi'\mid_{tA}$.\\

{\bf Remark 1.3}. Let $\beta:B\to X$ be an isomorphism and $x\in tX$. Then $nx=0$ for some positive integer $n$. Since $\beta$ is surjective, so there exist $b\in B$ such that $\beta(b)=x$. Hence, $\beta(nb)=0$. Since $\beta$ is injective, so $nb=0$. Therefore, $\beta\mid_{tB}:tB\to tX$ is an isomorphism.\\

Recall that two extensions $0 \to A \stk{\phi_{1}} B \stk{\psi_{1}} C \to 0$ and $0 \to A \stk{\phi_{2}} X \stk{\psi_{2}} C \to 0$ are said to be equivalent if there is an isomorphism $\beta:B\to X$ such that the following diagram
\[
\xymatrix{
0 \ar[r] & A \ar^{\phi_{1}}[r] \ar^{1_{A}}[d] & B \ar^{\psi_{1}}[r] \ar^{\beta}[d] & C \ar[r] \ar^{1_{C}}[d]
& 0 \\
0 \ar[r] & A \ar^{\phi_{2}}[r] &   X \ar^{\psi_{2}}[r] & C \ar[r] & 0
}
\]
is commutative.\\

{\bf Lemma 1.4}. An extension equivalent to a t-extension is a t-extension.\\

Proof. Let $$ E_1: 0 \to A \stk{\phi_1} B\stk{\psi_1} C \to 0 $$  and $$ E_2: 0 \to A \stk{\phi_2} X \stk{\psi_2} C \to 0 $$  be two equivalent extensions such that $E_1$ is a t-extension. Then, there is an isomorphism $\beta:B\to X$ such that $\beta\phi_1=\phi_2$ and $\psi_2\beta=\psi_1$. Let $x\in tC$. Since $E_1$ is t-extension, so $\psi_1(b)=x$ for some $b\in tB$. Hence, $\psi_2(\beta(b))=\psi_1(b)=x$. So, $\psi_2:tX\to tC$ is surjective. Now, let $\psi_2(x)=0$ for some $x\in tX$. By Remark 1.3, there exist $b\in tB$ such that $\beta(b)=x$. Hence, $\psi_1(b)=\psi_2(\beta(b))=0$. Since $E_1$ is t-extension, so $\phi_1(a)=b$ for some $a\in tA$. Consequently, $\phi_2(a)=\beta(\phi_1(a))=x$.\\

{\bf Remark 1.5}. Let $C$ and $A$ be two groups and $0\to A\stk{\phi_{1}} B_{1}\stk{\psi_{1}} C\to 0$ and  $0\to A\stk{\phi_{2}} B_{2}\stk{\psi_{2}} C\to 0$  two t-extensions of $A$ by $C$. An easy calculation shows that $0\to A\bigoplus A\stk{(\phi_{1}\bigoplus\phi_{2})} B_{1}\bigoplus B_{2}\stk{(\psi_{1}\bigoplus \psi_{2})} C\bigoplus C\to 0$ is a t-extension where $(\phi_{1}\bigoplus\phi_{2})(a_{1},a_{2})=(\phi_{1}(a_{1}),\phi_{2}(a_{2}))$ and $(\psi_{1}\bigoplus \psi_{2})(b_{1},b_{2})=(\psi_{1}(b_{1}),\psi_{2}(b_{2}))$. \\

{\bf Theorem 1.6}. Let $A$ and $C$ be two groups. Then, the class $Ext_{t}(C,A)$ of all equivalence classes of t-extensions of $A$ by $C$ is an subgroup of $Ext(C,A)$ with respect to the operation defined by $$[E_{1}]+[E_{2}]=[\nabla_{A}(E_{1}\bigoplus E_{2})\triangle_{C}]$$ where $E_{1}$ and $E_{2}$ are t-extensions of $A$ by $C$ and $\nabla_{A}$ and $\triangle_{C}$ are the diagonal and codiagonal homomorphism.\\

Proof. Clearly, $0\to A\to A\bigoplus C\to C\to 0$ is a t-extension. By Lemma 1.2 and Remark 1.5, $[E_{1}]+[E_{2}]\in Ext_{t}(C,A)$ for two t-extensions $E_{1}$ and $E_{2}$ of $A$ by $C$. So, $Ext_{t}(C,A)$ is a subgroup of $Ext(C,A)$.\\

{\bf Lemma 1.7}. {\it Let $A$ and  $C$ be two groups. Then, $Pext(C,A)\subseteq Ext_{t}(C,A)$.}\\

Proof. Let $0\to A\stk{\phi} B\stk{\psi} C\to 0$ be an element of $Pext(C,A)$. It is sufficient to show that $tB/t\phi(A)\cong t(B/\phi(A))$. Consider the map $\varphi:tB\to t(B/\phi(A))$ given by $b\longmapsto b+\phi(A)$. Clearly, $\varphi$ is a homomorphism. We show that $\varphi$ is surjective. Let $b+\phi(A)\in t(B/\phi(A))$. Then, there exists a positive integer $n$ such that $nb\in \phi(A)$. Since $\phi(A)$ is pure in $B$, so $nb=n\phi(a)$ for some $a\in A$. Hence, $n(b-\phi(a))=0$. This shows that $b-\phi(a)\in tB$ and $\varphi(b-\phi(a))=b+\phi(A)$. So $\varphi$ is surjective. We have $$Ker\varphi=\{b\in tB;b\in \phi(A)\}=\phi(A)\bigcap tB=t\phi(A)$$ Hence $tB/t\phi(A)\cong t(B/\phi(A))$.\\

{\bf Corollary 1.8}. {\it If $A$ is a divisible group or $C$ a torsion-free group, then $Pext(C,A)=Ext_{t}(C,A)=Ext(C,A)$.}\\

Proof. It is clear.\\

\section{Splitting of t-extensions}

In this section, we establish some conditions on $A$ and $C$ such that $Ext_{t}(C,A)=0$. We also determine the t-injective and t-projective groups in $\Re$.\\

{\bf Lemma 2.1}. {\it Let $A$ be a torsion-free group and $C$ a torsion group. Then, $Ext_{t}(C,A)=0$.}\\

Proof. Let $E:0\to A\stk{\phi} B\stk{\psi} C\to 0$ be a t-extension. Then $\psi_{|tB}:tB\to C$ is an isomorphism. Let $b\in B$. Then, $\psi(b)\in C$. So $\psi(b)=\psi(b')$ for some $b'\in tB$. Hence, $b-b'=\phi(a)$ for some $a\in A$. This follows that $B=\phi(A)+tB$. Since $\phi(A)$ is torsion-free and $tB$ torsion, so $\phi(A)\bigcap tB=0$. Hence, $B=\phi(A)\bigoplus tB$ and $E$ splits. \\

{\bf Corollary 2.2}. {\it  Let $A$ be a torsion-free group and $C$ a torsion group. Then, $Pext(C,A)=0$.}\\

Proof. It is clear by Lemma 1.7 and Lemma 2.1.\\

{\bf Lemma 2.3}. Let $A$ and $C$ be two torsion groups. Then $Ext(C,A)=Ext_{t}(C,A)$.\\

Proof. Let $A$ and $C$ be two torsion groups. It is clear that $Ext_{t}(C,A)\subseteq Ext(C,A)$. Suppose that $E:0\to A\to B\to C\to 0$ is an extension. Then, $B$ is a torsion group. Hence, $E$ is a t-extension.\\

{\bf Lemma 2.4}. {\it Let $C$ be a torsion group. Then, $Ext(C,tA)\cong Ext_{t}(C,A)$ for every group $A$.}\\

Proof. The exact sequence $0\to tA\stk{i}A\stk{\pi}A/tA\to 0$ induces the following exact sequence $$Hom(C,A/tA)\to Ext(C,tA)\stk{i_{\ast}} Ext(C,A)\stk{\pi_{\ast}}Ext(C,A/tA)\to 0$$ Note that $Hom(C,A/tA)=0$. By Lemma 2.1, $Ext_{t}(C,A/tA)=0$. Since $\pi_{\ast}(Ext_{t}(C,A))\subseteq Ext_{t}(C,A/tA)$, so $\pi_{\ast}(Ext_{t}(C,A))=0$. Hence, $Ext_{t}(C,A)\subseteq Ker\pi_{\ast}=i_{\ast}(Ext(C,tA))$. By Lemma 2.3, $Ext(C,tA)=Ext_{t}(C,tA)$. So, $i_{\ast}(Ext(C,tA))\subseteq Ext_{t}(C,A) $. Hence, $Ext(C,tA)\cong Ext_{t}(C,A)$. \\

{\bf Corollary 2.5}. {\it Let $A$ be a group. Then, $ Ext_{t}(Z(m),A)\cong tA/m(tA)$ for every positive integer $m$.\\

Proof. It is clear by Lemma 2.4 and [1,page 222].\\

{\bf Corollary 2.6}. {\it Let $C$ be a torsion group and $\{A_{i}:i\in I\}$ a collection of groups. If $I$ is finite, then $Ext_{t}(C,\Pi_{i\in I}A_{i})\cong \Pi_{i\in I}Ext_{t}(C,A_{i})$.\\

Proof. Let $I=\{1,...,n\}$ for some positive integer $n$. By Lemma 2.4, $Ext_{t}(C,\Pi_{i=1}^{n}A_{i})\cong Ext(C,t(\Pi_{i=1}^{n}A_{i}))$. Since $t(\Pi_{i=1}^{n}A_{i})=\Pi_{i=1}^{n}t(A_{i})$, so $Ext_{t}(C,\Pi_{i=1}^{n}A_{i})\cong \Pi_{i\in I}Ext(C,tA_{i})\cong \Pi_{i\in I}Ext_{t}(C,A_{i})$.\\

{\bf Remark 2.7}. In general, $Ext_{t}(C,\Pi_{i\in I}A_{i})\not \cong \Pi_{i\in I}Ext_{t}(C,A_{i})$. \\

{\bf Example 2.8}. Let $p$ be a prime and $H=\prod_{n=1}^{\infty}Z(p^{n})$. By Lemma 2.4, $Ext_{t}(Q/Z,H)\cong Ext(Q/Z,tH)$. Consider the following exact sequence $$(2.1)~~~~~0\to Ext(Q/Z,tH)\to Ext(Q/Z,H)\to Ext(Q/Z,H/tH)\to 0$$ By [1, Theorem 52.2] and Lemma 2.3, $$Ext(Q/Z,H)\cong \prod_{n=1}^{\infty}Ext(Q/Z,Z(p^{n}))\cong\prod_{n=1}^{\infty}Ext_{t}(Q/Z,Z(p^{n}))$$ If $Ext_{t}(Q/Z,H)\cong \prod_{n=1}^{\infty}Ext_{t}(Q/Z,Z(p^{n}))$,then $Ext(Q/Z,H)\cong Ext(Q/Z,tH)$. It follows from (2.1) that, $Ext(Q/Z,H/tH)=0$ which is a contradiction, since $H/tH$ is not divisible.\\

{\bf Lemma 2.9}. {\it Let $A$ be a torsion-free group. Then, $Ext(C/tC,A)\cong Ext_{t}(C,A)$ for every group $C$.}\\

Proof. The exact sequence $0\to tC\stk{i}C\stk{\pi}C/tC\to 0$ induces the following exact sequence $$Hom(tC,A)\to Ext(C/tC,A)\stk{\pi_{\ast}} Ext(C,A)\stk{i_{\ast}}Ext(tC,A)\to 0$$ Note that $Hom(tC,A)=0$. By Lemma 2.1, $Ext_{t}(tC,A)=0$. So, $i_{\ast}(Ext_{t}(C,A))\subseteq Ext_{t}(tC,A)=0$. Hence, $Ext_{t}(C,A)\subseteq Ker i_{\ast}=\pi_{\ast}(Ext(C/tC,A))$. By Corollary 1.8, $Ext(C/tC,A)=Ext_{t}(C/tC,A)$. So, $\pi_{\ast}(Ext(C/tC,A))\subseteq Ext_{t}(C,A) $. Hence, $Ext(C/tC,tA)\cong Ext_{t}(C,A)$.\\

{\bf Definition 2.10}. {\it Let $G$ be a group. We call $G$ a t-injective group in $\Re$ if for every t-extension $$ 0 \to A \stk{\phi} B \to C \to 0 $$ and a homomorphism $ f:A\to G$, there is a homomorphism $ \bar{f}:B\longrightarrow G$ such that $ \bar{f}\phi=f$.

We call $G$ a t-projective group in $\Re$ if for every t-extension $$ 0 \to A \to B \stk{\psi} C \to 0 $$ and a homomorphism $ f:G \to C$, there is a homomorphism $\bar{f}:G\to B$ such that $\psi\bar{f}=f$ .}\\

Recall that a group $A$ is said to be cotorsion if $Ext(Q,A)=0$.\\

{\bf Theorem 2.11}. {\it Let $A$ be a group. The following statements are equivalent:
\begin{enumerate}
\item $A$ is t-injective in $\Re$.
\item $Ext_{t}(C,A)=0$ for all $C\in \Re$.
\item $A\cong B\bigoplus D$ where $B$ is a torsion divisible group and $D$ a torsion-free cotorsion group.\\
\end{enumerate}

Proof. $1\Longrightarrow 2$: Let $A$ be a t-injective in $\Re$ and $ E: 0\to A \stk{\phi} B \stk{} C\to 0$ a t-extension of $A$ by $C$. Then, there is a homomorphism $\bar{\phi}:B \to G$ such that $\bar{\phi}\phi=1_G$. Consequently, $E$ splits.

$2\Longrightarrow 3$: Let $Ext_{t}(C,A)=0$ for every group $C$. So $Ext_{t}(Z(m),A)=0$ for every positive integer $m$. By Corollary 2.5, $m(tA)=tA$ for every positive integer $m$. So, $tA$ is divisible. Hence, $A\cong tA\bigoplus A/tA$. Therefore, $Ext(Q,A)\cong Ext(Q,tA)\bigoplus Ext(Q,A/tA)$. Since $Ext(Q,A)=Ext_{t}(Q,A)=0$, then $Ext(Q,A/tA)=0$. Hence, $A/tA$ is cotorsion. Now, we set $tA=B$ and $A/tA=D$.

$3\Longrightarrow 2$: Suppose that $A\cong B\bigoplus D$ where $B$ is a torsion divisible group and $D$ a torsion-free cotorsion group. Let $C$ be a group. Since $Ext(C,B)=0$, so $\pi_{2}^{*}:Ext(C,A)\to Ext(C,D)$ is an isomorphism. By Lemma 2.9, $Ext_{t}(C,D)\cong Ext(C/tC,D)$. Since $D$ is a cotorsion group, so $Ext(C/tC,D)=0$. Hence, $Ext_{t}(C,D)=0$. So, $\pi_{2}^{*}(Ext_{t}(C,A))\subseteq Ext_{t}(C,D)=0$. Since $\pi_{2}^{*}$ is an isomorphism, so $Ext_{t}(C,A)=0$.

$2\Longrightarrow 1$: Let $ E:0 \to G \stk{\phi} B \to C \to 0 $ be a t-extension and $f:G\to A$ a homomorphism. Then $f$ induces a pushout diagram
\[
\xymatrix{
E:0 \ar[r] & G \ar^{\phi}[r] \ar^{f}[d] & B \ar^{}[r] \ar^{}[d] & C \ar[r] \ar^{}[d]
& 0 \\
fE:0 \ar[r] & A \ar^{\mu}[r] &   (A\bigoplus B)/ H \ar^{}[r] & C \ar[r] & 0
}
\]

Where $H= \{(-f(a),\phi(a));a\in A\}$ and $ \mu :a \longmapsto (a,0)+H$. By Lemma 1.2, $fE$ is a t-extension and by assumption it splits. Hence $A$ is t-injective.\\

Recall that a group $A$ is called algebraically compact if and only if $Pext(X,A)=0$ for every group $X$ [1].\\

{\bf Corollary 2.12}. {\it A torsion-free, cotorsion group is algebraically compact.}\\

Proof. Let $A$ be a torsion-free, cotorsion group. By Theorem 2.11, $Ext_{t}(C,A)=0$ for every group $C$. Hence, $Pext(C,A)=0$ for every group $C$.\\

{\bf Theorem 2.13}. {\it Let $C$ be a group. The following statements are equivalent:
\begin{enumerate}
\item $C$ is t-projective in $\Re$.
\item $Ext_{t}(C,A)=0$ for all $A\in \Re$.
\item $C$ is a free group.\\
\end{enumerate}

Proof. $1\Longrightarrow 2$: Let $C$ be a t-projective in $\Re$ and $ E: 0\to A \stk{\phi} B \stk{} C\to 0$ a t-extension of $A$ by $C$. Then there is a homomorphism $\bar{\phi}:C \to B$ such that $\bar{\phi}\phi=1_C$. Consequently, $E$ splits.

$2\Longrightarrow 3$: Let $Ext_{t}(C,A)=0$ for every group $A$. By Lemma 1.7, $Pext(C,A)=0$ for every group $A$. Hence, $C$ is free group.

$3\Longrightarrow 1$: It is clear.\\

\section{Splitting of pure extensions in the category of locally compact abelian groups}

In this section, we determine the structure of a compact group $G$ such that $Pext(C,G)=0$ for all compact connected groups $C$.\\

{\bf Lemma 3.1}. {\it If $A$ is a compact totally disconnected group and $C$ a compact connected group, then $Pext(C,A)=0$.\\

Proof. By [7, Lemma 2.3], $Pext(C,A)\cong Pext(\hat{A},\hat{C})$. On the other hand, $\hat{A}$ and $\hat{C}$ are a discrete, torsion group and a discrete, torsion-free respectively [6, Theorem 24.25 and 24.26]. Hence, by Corollary 2.2,$Pext(\hat{A},\hat{C})=0$. So $Pext(C,A)=0$. \\

Let $G$ be a compact group. Then, $Ext(C,G)=0$ for every compact connected group $C$ if and only if $G\cong (R/Z)^{\sigma}$ [5, Corollary 3.4]. In the next Theorem, we show that $Pext(C,G)=0$ for every compact connected group $C$ if and only if $G\cong (R/Z)^{\sigma}\bigoplus H$ where $H$ is a compact totally disconnected group. \\

{\bf Theorem 3.2}. {\it Let $G$ be a compact group. Then, $Pext(C,G)=0$ for all compact connected groups $C$ if and only if $G\cong (R/Z)^{\sigma}\bigoplus H$ where $H$ is a compact totally disconnected group.\\

Proof. Let $G$ be a compact group and $C$ a compact connected group. Consider the exact sequence $0\to G_{0}\to G\to G/G_{0}$. By Proposition 4 of [2], we have the exact sequence $$Hom(C,G/G_{0})\to Ext(C,G_{0})\to Pext(C,G)\to Pext(C,G/G_{0})$$
 By Lemma 3.1, $Pext(C,G/G_{0})=0$. Since $G/G_{0}$ is totally disconnected, so $Hom(C,G/G_{0})=0$. It follows that $Pext(C,G)\cong Ext(C,G_{0})$. Now, Let $Pext(C,G)=0$ for all compact connected groups $C$. Then $Ext(C,G_{0})=0$ for all compact connected groups $C$. By [5, Corollary 3.4], $G_{0}\cong (R/Z)^{\sigma}$. So, the extension $0\to G_{0}\to G\to G/G_{0}\to 0$ splits. This shows that $G\cong(R/Z)^{\sigma}\bigoplus G/G_{0}$. The converse is clear.\\

{\bf Remark 3.3}. {\it Let $G$ be a compact group such that $Ext(C,G)=0$ for every compact connected group $C$. Then $Pext(C,G)=0$ for every compact connected group $C$. So by Theorem 3.2, $G\cong (R/Z)^{\sigma}\bigoplus H$ where $H$ is a compact totally disconnected group. Since $R/Z$ is a compact connected group, so $Ext(R/Z,H)=0$. Consider the following exact sequence $$Hom(R,H)\to Hom(Z,H)\to Ext(R/Z,H)=0$$ Since $R$ is a connected group and $H$ a totally disconnected group, so $Hom(R,H)=0$. Hence $H=0$ and $G\cong R/Z$.\\

\end{document}